\documentclass[11pt]{amsart}
\usepackage[hmargin=2.5cm,vmargin=2.5cm]{geometry}
\usepackage{amsfonts}
\usepackage{amsthm}
\usepackage[english]{algorithm2e}
\usepackage{amsmath}
\usepackage{amssymb}
\usepackage{hyperref}
\usepackage[T1]{fontenc}
\newtheorem{Theorem}{Theorem}[section]
\newtheorem{Definition}[Theorem]{Definition}
\newtheorem{Remark}[Theorem]{Remark}
\newtheorem{Example}[Theorem]{Example}
\newtheorem{Lemma}[Theorem]{Lemma}
\newtheorem{Proposition}[Theorem]{Proposition}
\newtheorem{Corollary}[Theorem]{Corollary}

\title{Structures and bimodules of simple Hom-alternative algebras}
\author[Sylvain Attan ]
       { Sylvain Attan}
\begin{document}
\maketitle
\begin{abstract} This paper is mainly devoted  to a structure study of   Hom-alternative algebras . Equivalent conditions for  Hom-alternative algebras being solvable, simple and semi-simple are displayed. Moreover some results about Hom-alternative bimodule  are found. 
\end{abstract}
\section{Introduction}
Hom-algebras are new classes of algebras which has been studied extensively in the litterature during the last decade. They are algebras where the identities defining the structure are twisted by a homomorphism and started from Hom-Lie algebras  \cite{HAR1}, \cite {dlsds1}, \cite{dlsds2}, \cite {dlsds3}, motivated by quasi-deformations of Lie algebras of vector fields, in particular q-deformations of Witt and Virasoro algebras. Hom-associative algebras were introduced in \cite{MAK3} while Hom-alternative and Hom-Jordan algebras are introduced in \cite{MAK1}\cite{YAU3}.

Questions on the structure of simple algebras in this or that variety
are one of the main questions in the theory of rings. This question for alternative algebras, has been studied by many authors. It turns out that the only simple alternative algebras which are not associative are $8$-dimensional algebras over their centers which are generalizations of the original algebra of Cayley numbers \cite{Schaf0}, \cite{mz1}, \cite{mz2}. Hence all semisimple alternative algebras are known. 
Similarly, it is very necessary to study simple Hom-algebras in Hom-algebras theory. In \cite{xc}, the authors gave the classification theorem about multiplicative simple Hom-Lie algebras. Inspired by this study, classification of mutiplicative simple Hom-Jordan algebras
 is obtained \cite{cy}.

Representations (or bimodules) and deformations are important tools in most parts of Mathematics and Physics. By means of the representation theory, we would be more aware of the corresponding algebras.
The study of bimodule  of Jordan algebras was initiated by N. Jacobson \cite{Jacob1}.  Subsequently the alternative case was considered by Schafer \cite{Schaf1}. Similarly, it is very necessary to study representation of Hom-algebras. Fortunately, representations of Hom-Lie algebras were introduced and studied in \cite{YS}, see also \cite{sba}. Later  the one of Hom-Jordan and Hom-alternative is presented \cite{sa} where some useful results are obtained. Moreover representations of simple Hom-Lie algebras \cite{xc}, \cite{xxl}  and the one of simple Hom-Jordan algebras \cite{cy} are introdued and studied in detail.
In this paper, basing on \cite{cy}  and on \cite{az} ,  we  will study structure and bimodule over simple Hom-alternative algebras.

The paper is organized as follows: In section two, we give the basics about Hom-alternative algebras and provide some new properties. Section three deals with simple and solvable Hom-alternative algebras. Here some  useful results are obtained for the next sections (see Lemma \ref{lkhi} and Theorem \ref{haat}). 
In section four we mainly prove relevant theorems which are about solvability, simpleness and semi-simpleness of Hom-alternative algebras 
(see theorem \ref{tes}, Theorem \ref{tcns1} and Theorem \ref{ts}).  
In section five, we prove Theorem \ref{bHaba} which is very important. It talks about the relationship between bimodules over Hom-alternative algebras of alternative type and the ones over their induced alternative algebras. Moreover, some relevant propositions about bimodule over Hom-alternative algebras also displayed as applications of Theorem \ref{bHaba}.

All vector spaces are assumed to be  over a fixed ground field $\mathbb{K}$ of characteristic $0.$
\section{Preliminaries}
We recall some basic notions, introduced in \cite{HAR1}, \cite{MAK3}, \cite{YAU4} related to Hom-algebras and while dealing of any binary operation  we will use juxtaposition in order to reduce the 
number of braces i.e., e.g., for $"\cdot ", $  $xy\cdot\alpha(z)$ means $(x\cdot y)\cdot\alpha(z).$ 
Also for the map $\mu: A^{\otimes 2}\longrightarrow A,$ we will write sometimes $\mu(a\otimes b)$ as $\mu(a,b)$ or $ab$  for $a,b\in A.$
\begin{Definition}
A Hom-module is a pair $(M,\alpha_M)$ consisting of a $\mathbb{K}$-module $M$ and 
a linear self-map $\alpha_M: M\longrightarrow M.$ A morphism 
$f: (M,\alpha_M)\longrightarrow (N,\alpha_N)$ of Hom-modules is a linear map 
 $f: M\longrightarrow N$ such that $f\circ\alpha_M=\alpha_N\circ f.$
\end{Definition}
\begin{Definition} (\cite{MAK3}, \cite{YAU4}) A Hom-algebra is a triple $(A,\mu,\alpha)$ in which $(A,\alpha)$ is a Hom-module, $\mu: A^{\otimes 2}\longrightarrow A$ is a linear map.
The Hom-algebra $(A,\mu,\alpha)$ is said to be  multiplicative if $\alpha\circ\mu=\mu\circ\alpha^{\otimes 2}$ (multiplicativity). A morphism 
$f: (A,\mu_A,\alpha_A)\longrightarrow (B,\mu_B,\alpha_B)$ of Hom-algebras is a morphism of the underlying Hom-modules such that $f\circ\mu_A=\mu_B\circ f^{\otimes 2}.$
\end{Definition}
In this paper, we will only consider multiplicative Hom-algebras.
\begin{Definition}\label{DefRef}
Let $(A,\mu,\alpha)$ be a Hom-algebra.\\
(1) The Hom-associator of $A$ is the linear map $as_A: A^{\otimes 3}\longrightarrow A$  defined as $as_A=\mu\circ(\mu\otimes\alpha-\alpha\otimes\mu).$ A multiplicative Hom-algebra $(A,\mu,\alpha)$ is said to be Hom-associative algebra if $as_A=0$\\
(2) A  Hom-alternative algebra \cite{MAK1} is a multiplicative Hom-algebra $(A,\mu,\alpha)$ that satisfies both
\begin{eqnarray}
as_A(x,x,y)&=&0  \mbox{ (left Hom-alternativity)}\label{lhi} \mbox { and }\\
as_A(x,y,y)&=&0  \mbox{ (right Hom-alternativity)}\label{rhi}
\end{eqnarray}
for all $x,y\in A.$\\
(3) Let $(A, \mu, \alpha)$ be a Hom-alternative algebra. A Hom-subalgebra of  $(A, \mu, \alpha)$ is a linear subspace $H$ of $A,$ which is closed for the multiplication $\mu$ and invariant by $\alpha,$ that is, $\mu(x,y)\in H$ and $\alpha(x)\in H$ for all $x,y\in H.$ If furthermore $\mu(a,b)\in H$ and $\mu(b,a)\in H$  for all $(a,b)\in A\times H,$ then $H$ is called a two-sided Hom-ideal of $A.$
\end{Definition}
\begin{Example} \label{ex1} The octonions algebra $\mathbb{O}$ also called Cayley Octaves or Cayley algebra is
8-dimensional  with  a basis $(e_0,e_1,e_2,e_3,e_4,e_5,e_6,e_7),$ where $e_0$ is the identity for the multiplication.
 This algebra is twisted into the eight-dimensional Hom-alternative algebra $\mathbb{O}_{\alpha}=(\mathbb{O}, \mu_1, \alpha)$ \cite{YAU3}  with the same basis $(e_0,e_1,e_2,e_3,e_4,e_5,e_6,e_7)$ where
$\alpha(e_0)=e_0,\ \alpha(e_1)=e_5,\ \alpha(e_2)=e_6,\ \alpha(e_3)=e_7,\ \alpha(e_4)=e_1,\alpha(e_5)=e_2,\ \alpha(e_6)=e_3,\ \alpha(e_7)=e_4$  and the multiplication table is:
$$
\begin{array}{|c|c|c|c|c|c|c|c|c|c}
\hline
\mu_1 &e_0&e_1&e_2&e_3&e_4&e_5&e_6&e_7\\
\hline
e_0 &e_0&e_5&e_6&e_7&e_1&e_2&e_3&e_4 \\
\hline
e_1 &e_5&-e_0&e_1&e_4&-e_6&e_3&-e_2&-e_7\\
\hline
e_2&e_6&-e_1&-e_0&e_2&e_5&-e_7&e_4&-e_3\\
\hline
e_3&e_7&-e_4&-e_2&-e_0&e_3&e_6&-e_1&e_5\\
\hline
e_4&e_1&e_6&-e_5&-e_3&-e_0&e_4&e_7&-e_2\\
\hline
e_5&e_2&-e_3&e_7&-e_6&-e_4&-e_0&e_5&e_1\\
\hline
e_6&e_3&e_2&-e_4&e_1&-e_7&-e_5&-e_0&e_6 \\
\hline
e_7&e_4&e_7&e_3&-e_5&e_2&-e_1&-e_6&-e_0\\
\hline
\end{array},\
$$
and  into the eight-dimensional Hom-alternative algebra  $\mathbb{O}_{\beta}=(\mathbb{O}, \mu_2, \beta)$ \cite{MAK1}  with the same basis $(e_0,e_1,e_2,e_3,e_4,e_5,e_6,e_7)$ where
$\beta(ei)=-ei$ for all $i\in\{0,\cdots,7\}$  and the multiplication table is:
$$
\begin{array}{|c|c|c|c|c|c|c|c|c|c}
\hline
\mu_2 &e_0&e_1&e_2&e_3&e_4&e_5&e_6&e_7\\
\hline
e_0 &e_0&-e_1&-e_2&-e_3&e_4&e_5&-e_6&e_7 \\
\hline
e_1 &-e_1&-e_0&e_4&e_7&e_2&-e_6&-e_5&e_3\\
\hline
e_2&-e_2&-e_4&-e_0&e_5&-e_1&e_3&e_7&e_6\\
\hline
e_3&-e_3&-e_7&-e_5&-e_0&-e_6&-e_2&-e_4&-e_1\\
\hline
e_4&e_4&-e_2&e_1&e_6&-e_0&e_7&-e_3&-e_5\\
\hline
e_5&e_5&e_6&-e_3&e_2&-e_7&-e_0&-e_1&e_4\\
\hline
e_6&-e_6&e_5&-e_7&e_4&e_3&e_1&-e_0&-e_2 \\
\hline
e_7&e_7&-e_3&-e_6&e_1&e_5&-e_4&e_2&-e_0\\
\hline
\end{array}
$$
Nor $\mathbb{O}_{\alpha},$ neither $\mathbb{O}_{\beta},$ are alternative algebras. Moreover, both $\alpha$ and $\beta$ are automorphisms of $\mathbb{O}.$ 
\end{Example}
Similarly as in \cite{MAK1}, it easy to prove the following:
\begin{Proposition} \label{chafaa} Let $(A,\mu, \alpha)$ be a Hom-alternative algebra and 
$\beta: A\longrightarrow A$ be a morphism of $(A,\mu,\alpha).$ Then $(A,\beta\circ\mu,\beta\circ\alpha)$ is a Hom-alternative algebra.  In particular, if $(A,\mu)$ is an alternative algebra and $\beta$ is a morphism of $(A,\mu),$ then $(A,\beta\circ\mu,\beta$ is a Hom-alternative algebra \cite{MAK1}.
\end{Proposition}
\begin{Definition}
Let $(A,\mu,\alpha)$ be a Hom-alternative algebra. If there is an alternative algebra $(A,\mu')$ such that $\mu=\alpha\circ\mu',$ we say that $(A,\mu,\alpha)$ is of alternative type and $(A,\mu')$ is its compatible alternative algebra or the untwist of $(A,\mu,\alpha).$
\end{Definition}
It is remarked \cite{MAK1} that a Hom-alternative algebra with an invertible twisting map has the compatible alternative algebra. More precisely, we get:
\begin{Corollary} \label{chaaa}
Let $(A,\mu,\alpha)$ be a Hom-alternative algebra where $\alpha$ is invertible then $(A,\mu'=\alpha^{-1}\circ\mu)$ is an alternative algebra and $\alpha$ is an automorphism with respect to $\mu'.$ Hence $(A,\mu,\alpha)$ is of alternative type and $(A,\mu'=\alpha^{-1}\circ\mu)$ is its compatible alternative algebra.
\end{Corollary}
{\bf Proof.} Let denote  by $as$ the associator of $(A,\mu'=\alpha^{-1}\circ\mu).$ We prove that $(A,\mu'=\alpha^{-1}\circ\mu)$ is an alternative algebra. Indeed,
\begin{eqnarray}
as(x,x,y)&=&\mu'(\mu'(x,x),y)-\mu'(x,\mu'(x,y))=\alpha^{-1}\circ\mu(\alpha^{-1}\circ\mu(x,x),y)-
\alpha^{-1}\circ\mu(x,\alpha^{-1}\circ\mu(x,y))\nonumber\\
&=&\alpha^{-1}\circ\mu(\alpha^{-1}\circ\mu(x,x),\alpha^{-1}\circ\alpha(y))-
\alpha^{-1}\circ\mu(\alpha^{-1}\circ\alpha(x),\alpha^{-1}\circ\mu(x,y))\nonumber\\
&=&\alpha^{-2}(\mu(\mu(x,x),\alpha(y))-\mu(\alpha(x),\mu(x,y)))=\alpha^{-2}\circ as_A(x,x,y)=\alpha^{-2}(0)=0.\nonumber
\end{eqnarray}
Similarly, we get $as(x,y,y)=0.$ Moreover, $\alpha$ is an automorphism with respect to $\mu'$. Indeed,
\begin{eqnarray}
\mu'(\alpha(x),\alpha(y))=\alpha^{-1}\circ\mu(\alpha(x),\alpha(y))=
\alpha^{-1}\circ\alpha\circ\mu(x,y)=\alpha\circ\mu'(x,y).\nonumber
\end{eqnarray}
\begin{Proposition}
Let $(A_1,\mu_1,\alpha_1)$ and $(A_2,\mu_2,\alpha_2)$ be two Hom-alternative algebras and $\varphi: $ be an invertible morphism of Hom-algebras. If $(A_1,\mu_1,\alpha_1)$ is of alternative type and $(A_1,\mu_1')$ is its compatible alternative algebra then $(A_2,\mu_2,\alpha_2)$ is of alternative type with compatible alternative algebra $(A_2,\mu_2'=
\varphi\circ\mu_1'\circ(\varphi^{-1}\otimes\varphi^{-1}))$ such that $\varphi: (A_1,\mu_1')\longrightarrow (A_2,\mu_2')$ is an algebra morphism.
\end{Proposition}
{\bf Proof.} We have $\alpha_2\circ\varphi=\varphi\circ\alpha_1$ and $\varphi$ defines $\mu_2$ by $\mu_2\circ\varphi^{\otimes 2}=\varphi\circ\mu_1$ since $\varphi$ is a morphism from $(A_1,\mu_1,\alpha_1)$ to $(A_2,\mu_2,\alpha_2).$ It is easy to show that to check that $(A_2,\mu_2')$ is an alternative algebra. Furthermore
\begin{eqnarray}
\mu_2\circ\varphi^{\otimes 2}=\varphi\circ\mu_1=\varphi\circ\alpha_1\circ\mu_1'
=\alpha_2\circ\varphi\circ\mu_1'=\alpha_2\circ\mu_2'\circ\varphi^{\otimes 2}\nonumber
\end{eqnarray}
We show that $(A_2,\mu_2')$ is an alternative algebra such that 
$\mu_2'=\varphi\circ\mu_1'\circ(\varphi^{-1}\otimes\varphi^{-1}).$ Let denote by $as_1$ and $as_2$ the associators of $(A_1,\mu_1)$ and $(A_2,\mu_2')$ respectively. Then 
\begin{eqnarray}
as_2(u,u,v)&=&\mu_2'(\mu_2'(u,u),v)-\mu_2'(u,\mu_2'(u,v))
=\varphi\circ\mu_1'(\varphi^{-1}\circ\varphi\circ\mu_1'(\varphi^{-1}(u),\varphi^{-1}(u)),\varphi^{-1}(v))\nonumber\\
&&-\varphi\circ\mu_1'(\varphi^{-1}(u),\varphi^{-1}\circ\varphi\circ\mu_1'(\varphi^{-1}(u),\varphi^{-1}(v)))=
\varphi\circ\mu_1'(\mu_1'(\varphi^{-1}(u),\varphi^{-1}(u)),\varphi^{-1}(v))\nonumber\\
&&-\varphi\circ\mu_1'(\varphi^{-1}(u),\mu_1'(\varphi^{-1}(u),\varphi^{-1}(v)))=\varphi\circ as_1(\varphi^{-1}(u),\varphi^{-1}(u),\varphi^{-1}(v))=\varphi(0)=0.\nonumber.
\end{eqnarray}
Similarly, we prove that $as_2(u,v,v)=0.$ Hence $(A,\mu_2')$ is an alternative algebra.\\

The following characterization was given for Hom-Lie algebra \cite{YS}
 and Hom-associative algebra \cite{az}.
 \begin{Proposition}
 Given two Hom-alternative algebras $(A,\mu_A,\alpha)$ and $(B,\mu_B,\beta),$ there is a Hom-alternative algebra $(A\oplus B, \mu_{A\oplus B}, \alpha+\beta),$ where the bilinear map $\mu_{A\oplus B}:(A\oplus B)^{\times 2}\longrightarrow (A\oplus B)$ is given by
 \begin{eqnarray}
 \mu_{A\oplus B}(a_1+b_1,a_2+b_2)=\mu_A(a_1,a_2)+\mu_B(b_1,b_2), \forall \ a_1,a_2\in A,\ \forall \ b_1,b_2\in B,\nonumber
 \end{eqnarray}
 and the linear  map $(\alpha+\beta): (A\oplus B)\longrightarrow (A\oplus B)$ is given by
 \begin{eqnarray}
 (\alpha+\beta)(a+b)&=& (\alpha(a)+\beta(b)) \ \forall \ (a,b)\in A\times B.
 \end{eqnarray}
 \end{Proposition}
 {\bf Proof.}
First, $(\alpha+\beta)$ is multiplicative with respect to $\mu_{A\oplus B}.$ Indeed,
\begin{eqnarray}
&&(\alpha+\beta)\circ(\mu_{A\oplus B})(a_1+b_1,a_2+b_2)=(\alpha+\beta)(\mu_A(a_1,a_2)+\mu_B(b_1,b_2))\nonumber\\
&&=\alpha\circ\mu_A(a_1,a_2)+\beta\circ\mu_B(b_1,b_2)
=\mu_A(\alpha(a_1),\alpha(a_2))+\mu_B(\beta(b_1),\beta(b_2)\nonumber\\
&&=\mu_{A\oplus B}(\alpha(a_1)+\beta(b_1),\alpha(a_2)+\beta(b_2))=
\mu_{A\oplus B}((\alpha+\beta)(a_1+b_1),(\alpha+\beta)(a_2+b_2))\nonumber
\end{eqnarray}
Secondly we prove the left Hom-alternativity (\ref{lhi}) for $A\oplus B$ as follows
\begin{eqnarray}
 &&as_{A\oplus B}(a_1+b_1,a_1+b_1, a_2+b_2)\nonumber\\
 &=&\mu_{A\oplus B}(\mu_{A\oplus B}(a_1+b_1,a_1+b_1),(\alpha+\beta)(a_2+b_2))\nonumber\\
&&-\mu_{A\oplus B}((\alpha+\beta)(a_1+b_1), \mu_{A\oplus B}(a_1+b_1,a_2+b_2))\nonumber\\
&&=\mu_{A\oplus B}(\mu_A(a_1,a_1)+\mu_B(b_1,b_1),(\alpha(a_2)+\beta(b_2))
\nonumber\\
&&-\mu_{A\oplus B}(\alpha(a_1)+\beta(b_1),\mu_A(a_1,a_2)+\mu_(b_1,b_2))
\nonumber\\
&=&\mu_A(\mu_A(a_1,a_1),\alpha(a_2))+\mu_B(\mu_B(b_1,b_1),\beta(b_2))\nonumber\\
&&-\mu_A(\alpha(a_1),\mu_A(a_1,a_2))-\mu_B(\beta(b_1),\mu_B(b_1,b_2))
\nonumber\\
&=&as_A(a_1,a_1,a_2)+as_B(b_1,b_1,b_2)=0\nonumber
\end{eqnarray}
Similarly, we prove the right Hom-alternativity (\ref{rhi}) for $A\oplus B.$ Hence $(A\oplus B, \mu_{A\oplus B}, \alpha+\beta)$ is a Hom-alternative algebra.
\begin{Proposition}
Let $(A,\mu_A,\alpha)$ and $(B,\mu_B,\beta)$  be two Hom-alternative algebras and $\varphi: A\rightarrow B$ be a linear map.  Denote 
by $\Gamma_{\varphi}\subset A\oplus B,$ the graph of $\varphi.$ Then 
$\varphi$ is a morphism from the Hom-alternative $(A,\mu_A,\alpha)$  to the Hom-alternative algebra $(B,\mu_B,\beta)$  if and only if its graph 
$\Gamma_{\varphi}$ is a Hom-subalgebra of $(A\oplus B, \mu_{A\oplus B}, \alpha+\beta).$
\end{Proposition}
{\bf Proof}
  Let $\phi: (A, \mu_A, \alpha)\longrightarrow (B, \mu_B ,\beta)$ be a morphism of Hom-alternative algebras. \\
Then we have for all $u, v\in A,$ 
$$ \mu_{A\oplus B}((u,\varphi(u)),(v,\varphi(v)))=(\mu_A(u,v),\mu_B(\varphi(u),\varphi(v)))=(\mu_A(u,v),\varphi(\mu_A(u,v))).$$ 
Thus the graph $\Gamma_{\varphi}$ is closed under the multiplication $\mu_{A\oplus B}.$
Furthermore since $\varphi\circ\alpha=\beta\circ\varphi$, we have $(\alpha\oplus\beta)(u, \varphi(u)) = (\alpha(u), 
\beta\circ\varphi(u)) = (\alpha(u), \varphi\circ\alpha(u)),$
which implies that $\Gamma_{\varphi}$ is closed under $\alpha\oplus\beta.$ Thus $\Gamma_{\varphi}$ is a Hom-subalgebra of 
$(A\oplus B, \mu_{A\oplus B} , \alpha\oplus \beta).$\\
Conversely, if the graph $\Gamma_{\varphi}\subset A\oplus B$ is a 
Hom-subalgebra of 
$(A\oplus B, \mu_{A\oplus B}, \alpha\oplus\beta),$ then we have
$$\mu_{A\oplus B}((u, \varphi(u)), (v, \varphi(v)))=(\mu_A(u, v) , \mu_B(\varphi(u), \varphi(v)) )\in\Gamma_{\varphi},$$ 
which implies that
$$\mu_B(\varphi(u), \varphi(v))=\varphi(\mu_A(u, v)).$$
Furthermore, $(\alpha\oplus\beta )(\Gamma_{\varphi})\subset\Gamma_{\varphi}$ implies
$$(\alpha\oplus\beta)(u, \varphi(u))=(\alpha(u), \beta\circ\varphi(u)) \in\Gamma_{\varphi} ,$$
which is equivalent to the condition $\beta\circ\phi(u)=\phi\circ\alpha(u),$ i.e. $ \beta\circ\varphi=\varphi\circ\alpha.$ Therefore, $\varphi$ is a
morphism of Hom-alternative algebras.\hfill $\square$
\section{Simple and solvable Hom-alternative algebras}
In this section, se study simple and solvable Hom-alternative algebras. This study is inspired by the study given in \cite{cy} and \cite{az}.
\begin{Definition}
Let $(A,\mu,\alpha)$ be a Hom-alternative algebra. Define its derived sequences as follow:
$$A^{(0)}=A,\ A^{(1)}=\mu(A,A),\ A^{(2)}=\mu(A^{(1)},A^{(1)}), \cdots, A^{(k)}=\mu(A^{(k-1)},A^{(k-1)}),\cdots$$
\end{Definition}
The following elementary result will be very useful.
\begin{Lemma} \label{lemma1}
Let $(A,\mu,\alpha)$ be a Hom-alternative algebra and $k\in\mathbb{N}.$ Then 
$$\cdots\subseteq A^{(k+1)}\subseteq A^{(k)}\subseteq A^{(k-1)}\subseteq\cdots \subseteq A^{(2)}\subseteq A^{(1)}\subseteq A$$
and $A^{(k)}$ is a two-sided Hom-ideal of $(A,\mu,\alpha).$
\end{Lemma}
{\bf Proof.}  First, it is clear that $A^{(1)}\subseteq A^{(0)}=A.$ Next let $k\in\mathbb{N}$ and assume that $A^{(k)}\subseteq A^{(k-1)}.$ Then, we have $$A^{(k+1)}=\mu(A^{(k)},A^{(k)})\subseteq \mu(A^{(k-1)},A^{(k-1)})=A^{(k)}$$
 To prove that $A^{(k)}$ is a two-sided Hom-ideal of $(A,\mu,\alpha)$ for every $k\in\mathbb{N},$ by the inclusion condition, it suffices to prove the case $n=1.$ If note that $A^{(1)}\subseteq A,$ we have first
 $$\alpha(A^{(1)})=\alpha(\mu(A,A))=\mu(\alpha(A),\alpha(A))\subseteq 
 \mu(A,A)=A^{(1)}$$  and next
 $$\mu(A^{(1)},A)\subseteq\mu(A,A)=A^{(1)} \mbox{ and } \mu(A,A^{(1)})\subseteq\mu(A,A)=A^{(1)}.$$
Thus $A^{(1)}$ is a two-sided Hom-ideal of $(A,\mu,\alpha).$ 
\begin{Definition}
Let $(A,\mu,\alpha)$  be a  Hom-alternative algebra. Then  $(A,\mu,\alpha)$ is said to be solvable if 
there exists $n\in\mathbb{N}^*$ such that $A^{(n)}=\{0\}.$
\end{Definition}
We get the following example of solvable and non solvable Hom-alternative algebras respectively.
\begin{Example}
\begin{enumerate}
\item Consider the $3$-dimensional Hom-alternative algebra 
$(A,\mu,\alpha)$ with basis $(e_1,e_2,e_3)$ where 
$\mu(e_2,e_2)=e_1,\ \mu(e_2,e_3)=e_1,\  \mu(e_3,e_2)=e_1,\ \mu(e_3,e_3)=e_1,\ \mu(e_3,e_2)=e_1$ and $\alpha(e_1)=e_1,\ \alpha(e_2)=e_1+e_2,\ \alpha(e_3)=e_3.$ Actually, $(A,\mu,\alpha)$ is a 
Hom-associative algebra ( see \cite{az}, Theorem 3.6, Hom-algebra ${A}_7^3$ ). Then $A^{(1)}$ is an one-dimensional Hom-alternative algebra 
generated by $(e_1)$ defined as follows: $\mu(e_1,e_1)=0$ and $\alpha(e_1)=e_1.$ It follows that $A^{(2)}=\{0\}$ and therefore 
$(A,\mu,\alpha)$ is solvable.
\item Consider  in the Example \ref{ex1}, the Hom-alternative algebras $\mathbb{O}_{\alpha}$ and $\mathbb{O}_{\beta}.$ From their multiplication tables, we get $(\mathbb{O}_{\alpha})^{(1)}=\mathbb{O}_{\alpha}$ and 
$(\mathbb{O}_{\beta})^{(1)}=\mathbb{O}_{\beta}.$ Hence for every $k\in\mathbb{N}$  we have $(\mathbb{O}_{\alpha})^{(k)}=\mathbb{O}_{\alpha}$ and $(\mathbb{O}_{\beta})^{(k)}=\mathbb{O}_{\beta}.$ It follows that nor 
 $\mathbb{O}_{\alpha}$ and neither $\mathbb{O}_{\beta}$ are solvable.
\end{enumerate}
\end{Example}
\begin{Definition}
Let $(A,\mu,\alpha)$ ($\alpha\neq 0$) be a non trivial Hom-alternative algebra.
\begin{enumerate}
\item $(A,\mu,\alpha)$ is said to be a simple Hom-alternative algebra if $A^{(1)}\neq\{0\}$ and it has no proper two-sided Hom-ideal.
\item $(A,\mu,\alpha)$ is said to be a semi-simple Hom-alternative algebra if $$A=A_1\oplus A_2\oplus A_2\oplus\cdots\oplus A_p$$
where $A_i$ ($1\leq i\leq p$) are simple two-sided Hom-ideals of $(A,\mu,\alpha).$
\end{enumerate}
\end{Definition}
Let give the following example of  non simple Hom-alternative algebra.
\begin{Example}
Consider the three-dimensional Hom-alternative algebra 
$(A, \mu,\alpha)$ over $\mathbb{K}$ with basis $(e_1,e_2,e_3)$ defined by $\mu(e_1,e_1)=e_1,\ \mu(e_2,e_2)=e_2,\ \mu(e_3,e_3)=e_1,\  \mu(e_1,e_3)=
\mu(e_3,e_1)=-e_3$ and $\alpha(e_1)=e_1,\  \alpha(e_3)=-e_3.$ Actually, $(A,\mu,\alpha)$ is a Hom-associative algebra ( see \cite{az}, Theorem 3.12, Hom-algebra ${A'}_3^3$ ). Consider the subspace $I=span(e_1,e_3)$ of $A.$ Then one can observe that $(I,\mu,\alpha)$ is a proper two-sided Hom-ideal of $(A, \mu,\alpha).$ Hence the Hom-alternative $(A, \mu,\alpha)$ is not simple.
\end{Example}
We have the following elementary result which will be used in next sections.
\begin{Proposition}\label{proposition1}
Let $(A,\mu,\alpha)$  be a simple Hom-alternative algebra. Then for each
$k\in\mathbb{N},$ $A^{(k)}=A.$
\end{Proposition}
{\bf Proof.} Thank to the definition of $A^{(k)},$ it suffices to prove that $A^{(1)}=A.$\\
By the simpleness of $(A,\mu,\alpha),$ we have $A^{(1)}\neq\{0\}.$ Moreover by Lemma \ref{lemma1}, $A^{(1)}$ is a two-sided Hom-ideal of $(A,\mu,\alpha)$ which has no proper two-sided Hom-ideal, then $A^{(1)}=A.$\\

The following lemma is useful for next results.
\begin{Lemma} \label{lkhi} Let $(A,\mu,\alpha)$ be a Hom-alternative algebra. Then 
$(Ker(\alpha),\mu,\alpha)$ is a two-sided Hom-ideal of $(A,\mu,\alpha),$.
\end{Lemma}
{\bf Proof.} Obvious, $\alpha(x)=0\in Ker(\alpha)$ for all $x\in Ker(\alpha).$ Next, let $x,z \in A$ and $y\in Ker(\alpha).$ Then 
$\alpha(\mu(x,y))=\mu(\alpha(x),\alpha(y))=\mu(\alpha(x),0)=0$ and 
$\alpha(\mu(y,z))=\mu(\alpha(y),\alpha(z))=\mu(0,\alpha(z))=0.$ Thus 
$\mu(x,y)\in Ker(\alpha)$ and $\mu(y,z)\in Ker(\alpha)$ and it follows that  $(Ker(\alpha),\mu,\alpha)$ is a two-sided Hom-ideal.
\begin{Proposition}\label{haat}
Let $(A,\mu,\alpha)$ be a finite dimensional simple Hom-alternative algebra. Then  the Hom-alternative algebra is of alternative type and 
$\alpha$ is an automorphism of both $(A,\mu,\alpha)$ and its induced algebra.
\end{Proposition}
{\bf Prof.} By Lemma \ref{lkhi}, $Ker(\alpha)$ is a two-sided Hom-ideal of the simple Hom-alternative algebra $(A,\mu,\alpha).$ Therefore $Ker(\alpha)=\{0\}$ or 
$Ker(\alpha)=A.$ Since the Hom-alternative algebra is non trivial, it follows that $Ker(\alpha)=\{0\}$ and $\alpha$ is an automorphism. 
Thus, $A$ is of alternative type( see corollary \ref{chaaa}).\\
Let $(A,\mu'=\alpha^{-1}\circ\mu)$ be the induced (the compatible) alternative algebra of the simple Hom-alternative algebra $(A,\mu,\alpha).$ We have $$\alpha\circ\mu'=\alpha\circ\alpha^{-1}\circ\mu=
\alpha^{-1}\circ\mu\circ\alpha^{\otimes 2}=\mu'\circ\alpha^{\otimes 2}$$    i.e. $\alpha$ is both automorphism of $(A,\mu')$ and $(A,\mu,\alpha).$\\

As Hom-associative algebras case \cite{az}, by the above proposition, there exists an induced alternative algebra of any simple Hom-alternative algebra $(A,\mu,\alpha)$ and $\alpha$ is an automorphism of the induced alternative algebra. Moreover, their products are mutually  determined.
\begin{Theorem} \label{fisHa}
Two finite dimensional simple Hom-alternative algebras $(A_1,\mu_1,\alpha)$ and $(A_2,\mu_2,\beta)$ are isomorphic if and only if there exists an alternative algebra 
isomorphism $\varphi: A_1\rightarrow A_2$ (between their induced alternative algebras) which renders conjugate the two alternative algebra automorphisms $\alpha$ and $\beta$ that is
$\varphi\circ\alpha=\beta\circ\varphi.$ 
\end{Theorem}
{\bf Proof.} Since $(A_1,\mu_1,\alpha)$ and $(A_2,\mu_2,
\beta)$ are finite dimensional simple Hom-alternative algebras, they are of alternative-type. Let $(A_1,\mu_1')$ and $(A_2,\mu_2')$ be their induced alternative algebras respectively. Suppose that $\varphi: (A_1,\mu_1,\alpha)\rightarrow(A_2,\mu_2, \beta)$ is an isomorphism of Hom-alternative algebras, then $\varphi\circ\alpha=\beta\circ\varphi,$ thus 
$\beta^{-1}\circ\varphi=\varphi\circ\alpha^{-1}.$ Moreover,
$$\varphi\circ\mu_1'=\varphi\circ\alpha^{-1}\circ\alpha\circ\mu_1'=
\varphi\circ\alpha^{-1}\circ\mu_1=\beta^{-1}\circ\varphi\circ\mu_1=
\beta^{-1}\circ\mu_2\circ\varphi^{\otimes 2}=\mu_2'\circ\varphi^{\otimes 2}.$$
So, $\varphi$ is an isomorphism between the induced Hom-alternative algebras.

On the other hand, if there exists an  
isomorphism between the induced alternative algebras  satisfying
$\varphi\circ\alpha=\beta\circ\varphi,$ then
$$\varphi\circ\mu_1=\varphi\circ\alpha\circ\mu_1'=
\beta\circ\varphi\circ\mu_1'=\beta\circ\mu_2'\circ\varphi^{\otimes 2}=
\mu_2\circ\varphi^{\otimes 2}.$$
\section{Structures of Hom-alternative algebras} 
In this section, we discuss about necessary and sufficient conditions for 
Hom-alternative algebras to be solvable, simple and semi-simple. This study is inspired by the one done in \cite{cy}.  
In the classical case, we already know that every simple alternative
algebra is either an associative algebra or a Cayley-Dickson algebra over its center \cite{Schaf0}, \cite{mz1}. As it turns out, there is one nonassociative simple alternative algebras. Recall that the Cayley-Dickson algebras are a sequence $A_0,\ A_1,\ \cdots$ of non-associative $\mathbb{R}$-algebras with involution. The first few are familiar: $A_0=\mathbb{R},\ A_1=\mathbb{C},\ A_2=\mathbb{H}$\ (the quaternions) and $A 3=\mathbb{O}$\ (the octonions). Each algebra $A_n$ is constructed from the previous one $A_{n-1}$ by a doubling procedure . 
The first three Cayley-Dickson algebras  are associative and  and $\mathbb{O}$ is the only one non-associative alternative Caylay-Dickson algebra.  Alternativity fails in the higher Cayley-Dickson algebras. Basing on this fact, we will give in this section an example of non Hom-associative simple Hom-alternative algebra. First let recall the following:
\begin{Proposition} \cite{sa}\label{haq}.
Let $(A,\mu,\alpha)$ be a Hom-alternative algebra and $I$ two-sided Hom-ideal of $(A,\mu,\alpha).$ Then $(A/I,\bar{\mu},\bar{\alpha})$ is a Hom-alternative algebra  where $\bar{\mu}(\bar{x},\bar{y})=\overline{\mu(x,y)}$ and $\bar{\alpha}(\bar{x})=\bar{\alpha(x)}$
for all $\bar{x},\bar{y}\in A/I.$
\end{Proposition}
{\bf Proof.}  First, note that the multiplicativity of $\bar{\mu}$ with respect to $\bar{\alpha}$ follows from the one of $\mu$ with respect to $\alpha.$ Next, pick $\bar{x},\bar{y}\in A/I.$ Then the left alternativity 
 (\ref{lhi} in $(A/I,\bar{\mu},\bar{\alpha})$ is proved as follows
 \begin{eqnarray}
as_{A/I}(\bar{x},\bar{x},\bar{y})&=&\bar{u}(\bar{u}(\bar{x},\bar{x}),\bar{\alpha}(\bar{y}))-\bar{\mu}(\bar{\alpha}(\bar{x}),\bar{\mu}(\bar{x},
\bar{y})\nonumber\\
&=&\overline{\mu(\mu(x,x)\alpha(y))-\mu(\alpha(x),\mu(x,y))}=
\overline{as_A(x,x,y))}=\bar{0}.\nonumber
 \end{eqnarray}
 Hence we get (\ref{lhi}) for $(A/I,\bar{\mu},\bar{\alpha}).$ Similarly, we get (\ref{rhi}) and therefore $(A/I,\bar{\mu},\bar{\alpha})$ is a Hom-alternative algebra.
\begin{Corollary}
Let $(A,\mu,\alpha)$ be a finite dimensional Hom-alternative algebra such that $\alpha^2=\alpha.$ Then the Hom-alternative algebra $(A/Ker(\alpha),\bar{\mu},\bar{\alpha})$ is of alternative type.
\end{Corollary}
{\bf Proof.} It is clear that $(A/Ker(\alpha),\bar{\mu},\bar{\alpha})$ is a Hom-alternative algebra by Proposition \ref{haq} since by Lemma \ref{lkhi}, $Ker(\alpha)$ is a two-sided Hom-ideal of $A.$\\
$\bullet$ If $\alpha$  is invertible i.e. $Ker(\alpha)=\{0\},$ then 
$(A/Ker(\alpha),\bar{\mu},\bar{\alpha})=(A,\mu,\alpha)$ and $(A/Ker(\alpha),\bar{\mu},\bar{\alpha})$ is a Hom-alternative algebra 
 of alternative type (see Corollary \ref{chaaa}).\\
$\bullet$ If $\alpha$ is not invertible, then $Ker(\alpha)\neq\{0\}.$ Therefore we have to show that $\bar{\alpha}$ is invertible on the Hom-alternative algebra $(A/Ker(\alpha),\bar{\mu},\bar{\alpha}).$ Assume that $\bar{x}\in Ker(\bar{\alpha}).$ Then $\bar{\alpha(x)}=\bar{\alpha}(\bar{x})=\bar{0},$ i.e. $\alpha(x)\in Ker(\alpha).$ Since $\alpha^2=\alpha,$ we have 
$$\alpha(x)=\alpha^2(x)=\alpha(\alpha(x))=0,$$
which means that $x\in Ker(\alpha),$ i.e. $\bar{x}=\bar{0}.$ If follows that $\bar{\alpha}$ is invertible and thank to the Corollary \ref{chaaa}, 
the Hom-alternative algebra $(A/Ker(\alpha),\bar{\mu},\bar{\alpha})$ is
 of alternative type.
 \begin{Theorem}\label{tes}
 Let $(A,\mu,\alpha)$ be a  Hom-alternative algebra such that $\alpha$ 
 is invertible. Then $(A,\mu,\alpha)$  is solvable if and only if its induced alternative algebra $(A,\mu')$  is solvable.
 \end{Theorem}
 {\bf Proof.} Let $(A,\mu,\alpha)$ be a Hom-alternative algebra such that $\alpha$ is invertible. Denote the derived sequences of $(A,\mu')$ and 
 $(A,\mu,\alpha)$ by 
 $A^{(i)}$ and $A_{\alpha}^{(i)}$ ($i=1,2,\cdots$) respectively.\\
 Suppose that $(A,\mu')$ is solvable. Then there exists $p\in\mathbb{N}^*$ such that 
 $A^{(p)}=\{0\}.$ Note that
 \begin{eqnarray}
 A_{\alpha}^{(1)}=\mu(A,A)=\alpha\circ\mu'(A,A)=\alpha(A^{(1)}),\nonumber\\
 A_{\alpha}^{(2)}=\mu( A_{\alpha}^{(1)}, A_{\alpha}^{(1)})=\mu (\alpha(A^{(1)}), \alpha(A^{(1)}))=\alpha^2\circ\mu'( A^{(1)}, A^{(1)})=\alpha^2(A^{(2)}),\nonumber
 \end{eqnarray}
 so by induction  $A_{\alpha}^{(p)}=\alpha^p(A^{(p)}).$ It follows that 
 $A_{\alpha}^{(p)}=\{0\},$ which means that $(A,\mu,\alpha)$ is solvable.
 
 On the other hand, assume that $(A,\mu,\alpha)$ is solvable. Then there exists $q\in\mathbb{N}^*$ such that $A_{\alpha}^{(q)}=\{0\}.$ By the above proof, we get $A_{\alpha}^{(q)}=\alpha^q(A^{(q)}).$ Note that $\alpha^q$ is invertible since $\alpha$ is, then $A^{(q)}=\{0\}$ that is $(A,\mu')$ is solvable.
 \begin{Lemma}\cite{cy}\label{lac1} Let $\mathcal{A}$ be an algebra over a field 
 $\mathbb{K}$ that has the unique decomposition of direct sum of simple ideals $\mathcal{A}=\oplus_{i=1}^s \mathcal{A}_i$ where the $\mathcal{A}_i$ are not isomorphic to each other and $\alpha\in Aut(\mathcal{A}).$ Then $\alpha(\mathcal{A}_i)=\mathcal{A}_i$ ($i=1,2,\cdots, s$).
 \end{Lemma}
 \begin{Theorem} \label{tcns1}
 \begin{enumerate}
 \item Let $(A,\mu,\alpha)$ be a finite dimensional simple Hom-alternative algebra. Then its induced alternative algebra $(A,\mu')$ is semi-simple. Moreover, $(A,\mu')$ can be decomposed into direct sum of isomorphic simple ideals. In addition, $\alpha$ acts simply transitively on simple ideals of the induced alternative algebra.
 \item Let $(A,\mu')$ be a simple alternative algebra and $\alpha\in Aut(A).$ Then $(A,\mu=\alpha\circ\mu',\alpha)$ is a simple Hom-alternative algebra.
 \end{enumerate}
 \end{Theorem}
 {\bf Proof.} 
 \begin{enumerate}
 \item Thank to Corollary \ref{chaaa}, $\alpha$ is both automorphism 
  with respect to $\mu'$ and $\mu.$\\
  Assume that $A_1$ is the maximal solvable two-sided ideal of $(A,\mu').$
  Then there exists $p\in\mathbb{N}^*$ such that $A_1^{(p)}=\{0\}.$ Since
   $$\mu'(A,\alpha(A_1))=\mu'(\alpha(A),\alpha(A_1))=\alpha(\mu'(A,A_1))
  \subseteq \alpha(A_1),$$
  $$\mu'(\alpha(A_1),A)=\mu'(\alpha(A_1),\alpha(A))=\alpha(\mu'(A_1,A))
  \subseteq \alpha(A_1)\mbox{ and } (\alpha(A_1))^{(p)}=\alpha(A_1^{(p)})=\{0\},$$
  we obtain that $\alpha(A_1)$ is also a solvable two-sided ideal of $(A,\mu').$ Then $\alpha(A_1)\subseteq A_1.$ Moreover
  $$\mu(A,A_1)=\alpha(\mu'(A,A_1))\subseteq \alpha(A_1)\subseteq A_1
  \mbox{ and } \mu(A_1,A)=\alpha(\mu'(A_1,A))\subseteq \alpha(A_1)\subseteq A_1.$$ It follows that 
  $A_1$ is a two-sided Hom-ideal of $(A,\mu,\alpha)$ and then 
  $A_1=\{0\}$ or $A_1=A$ since $(A,\mu,\alpha)$ is simple.
  
  If $A_1=A$, thank to the proof of Theorem \ref{tes}, we obtain
  $$A_{\alpha}^{(p)}=\alpha^p(A^{(p)})=\alpha^p(A_1^{(p)})=\{0\}.$$ 
  On the other hand, by the simpleness of $(A,\mu,\alpha),$ we have by Proposition \ref{proposition1}, $A_{\alpha}^{(p)}=A,$ which is a contradiction. It follows that $A_1=\{0\}$ and thus, $(A,\mu')$ is semi-simple. 
  
  Now, by the semi-simpleness of $(A,\mu'),$ we have $A=\oplus_{i=1}^s A_i$ where for all 
  $i\in\{1,\cdots, s\},$ $A_i$ is a simple two-sided ideal of $(A,\mu').$ Since there may be isomorphic alternative algebras among $A_1,\cdots, A_s,$ we can rewrite $A$ as follows:
  $$A=A_{11}\oplus A_{12}\oplus\cdots\oplus A_{1m_1}\oplus A_{21}\oplus A_{22}\oplus\cdots\oplus A_{2m_2}\oplus\cdots\oplus A_{t1}\oplus A_{t2}\oplus\cdots A_{tm_t},$$ where $(A_{ij},\mu')\cong (A_{ik},\mu'),$ $1\leq j,k\leq m_i,$ $i=1,2,\cdots,t.$ Thank to Lemma \ref{lac1}, we have 
  $$\alpha(A_{i1}\oplus A_{i2}\oplus\cdots\oplus A_{im_i})=A_{i1}\oplus A_{i2}\oplus\cdots\oplus A_{im_i}.\  \mbox{Then }$$ $$\mu(A_{i1}\oplus A_{i2}\oplus\cdots\oplus A_{im_i},A)=\alpha(\mu'(A_{i1}\oplus A_{i2}\oplus\cdots\oplus A_{im_i},A))\subseteq$$ 
  $$\alpha(A_{i1}\oplus A_{i2}\oplus\cdots\oplus A_{im_i})=A_{i1}\oplus A_{i2}\oplus\cdots\oplus A_{im_i}\ \mbox{ and }$$ 
  
  $$\mu(A,A_{i1}\oplus A_{i2}\oplus\cdots\oplus A_{im_i})=\alpha(\mu'(A, A_{i1}\oplus A_{i2}\oplus\cdots\oplus A_{im_i}))\subseteq$$ 
  $$\alpha(A_{i1}\oplus A_{i2}\oplus\cdots\oplus A_{im_i})=A_{i1}\oplus A_{i2}\oplus\cdots\oplus A_{im_i}.$$
  It follows that $A_{i1}\oplus A_{i2}\oplus\cdots\oplus A_{im_i}$ are two-sided Hom-ideals of $(A,\mu,\alpha)$. Then the simpleness of $(A,\mu,\alpha)$ implies $A_{i1}\oplus A_{i2}\oplus\cdots\oplus A_{im_i}=\{0\}$ or $A_{i1}\oplus A_{i2}\oplus\cdots\oplus A_{im_i}=A.$ Therefore, all but one of $A_{i1}\oplus A_{i2}\oplus\cdots\oplus A_{im_i}$ must be equal to $A.$ Without lost of generality, assume that 
  $$A=A_{11}\oplus A_{12}\oplus\cdots\oplus A_{1m_1}.$$
   If $m_1=1$ then $(A,\mu')$ is simple. Else, 
   $$\alpha(A_{1p})=A_{1l}\ (1\leq l\neq p\leq m_1)$$ since, if 
  $$\alpha(A_{1p})=A_{1p}\ (1\leq p\leq m_1),$$
   then $A_{1p}$ would be a non trivial two-sided Hom-ideal of $(A,\mu,\alpha)$ which contradicts the simpleness of $(A,\mu,\alpha).$
   
In addition, it is clear that  $A_{11}\oplus\alpha(A_{11})\oplus\alpha^2(A_{11})\oplus\cdots
  \alpha^{m_1-1}(A_{11})$ is a two-sided Hom-ideal of $(A,\mu,\alpha).$ Therefore 
  $$A=A_{11}\oplus\alpha(A_{11})\oplus\alpha^2(A_{11})\oplus\cdots
  \alpha^{m_1-1}(A_{11}).$$
   In other words, $\alpha$ acts simply transitively on simple ideals of the induced alternative algebra.
   \item By Proposition \ref{chafaa}, it is clear that $(A,\mu,\alpha)$ is a Hom-alternative algebra.
   
   Assume that $A_1$ is a non-trivial two-sided Hom-ideal of $(A,\mu,\alpha),$ then we get 
   $$\mu'(A,A_1)=\alpha^{-1}(\mu(A,A_1))\subseteq\alpha^{-1}(A_1)\subseteq A_1 \mbox{ and } \mu'(A_1,A)=\alpha^{-1}(\mu(A_1,A))\subseteq\alpha^{-1}(A_1)\subseteq A_1.$$ Its follows that $A_1$ is a non trivial two-sided ideal of $(A,\mu'),$ contradiction. 
   It follows that $(A,\mu,\alpha)$ has no trivial ideals.\\
   If $\mu(A,A)=\{0\},$ then, 
   $$\mu'(A,A)=\alpha^{-1}(\mu(A,A))=\{0\},$$
   which is in contradicts the fact that $(A,\mu')$ is simple.
   It follows that $(A,\mu,\alpha)$ is simple.
 \end{enumerate}
 Next, we will give an example of non Hom-associative simple Hom-alternative algebra using some results about Cayley-Dickson algebras.
 \begin{Example} Consider Hom-alternative algebras $\mathbb{O}_{\alpha}=(\mathbb{O},\mu_1,\alpha)$ and $\mathbb{O}_{\beta}=(\mathbb{O},\mu_2,\beta)$ (see Example \ref{ex1}) obtained from the octonions algebra $\mathbb{O}$ which is an eight dimensional simple non-associative alternative algebra. Since $\alpha\in Aut(\mathbb{O})$  and $\beta\in Aut(\mathbb{O}),$ thank to Theorem \ref{tcns1} (2), both $\mathbb{O}_{\alpha}$ and $\mathbb{O}_{\beta}$ are eight dimensional simple Hom-alternative algebras.
 \end{Example}
 \begin{Proposition} \label{psha}
 The eight dimensional simple Hom-alternative algebras $\mathbb{O}_{\alpha}$ and $\mathbb{O}_{\beta}$ are not isomorphic
 \end{Proposition}
 {\bf Proof.} Suppose that the simple Hom-alternative algebras $\mathbb{O}_{\alpha}$ and $\mathbb{O}_{\beta}$ are  isomorphic. Then by Theorem \ref{fisHa},  there exists an alternative algebra 
isomorphism $\varphi: \mathbb{O}\rightarrow \mathbb{O}$ (between their induced alternative algebras) such that $\varphi\circ\alpha=\beta\circ\varphi.$ This means by the definition  of $\beta$ that for all $i\in\{0,\cdots,7\},$ 
$\varphi(\alpha(e_i))=-\varphi(e_i)$ i.e. $\alpha=Id_{\mathbb{O}}$ (Contradiction).
\begin{Remark} We know  in the classical case that there exits only one example of a nonassociative simple alternative algebra that is, the Cayley-Dickson algebra over its center \cite{Schaf0}, \cite{mz1}. In the Hom-algebra setting, Proposition \ref{psha}, clearly proves that there is more than one simple non Hom-associative Hom-alternative algebras.
\end{Remark}
\begin{Theorem} \label{ts}
 \begin{enumerate}
 \item Let $(A,\mu,\alpha)$ be a finite dimensional semi-simple Hom-alternative algebra. Then   
  $(A,\mu,\alpha)$ is of alternative type and its induced alternative algebra $(A,\mu')$ is also semi-simple.
 \item Let $(A,\mu')$ be a semi-simple alternative algebra such that $A$ has a decomposition $A=\oplus_i^s A_i$  where $A_i$($1\leq i\leq s$) are simple two-sided ideal of $(A,\mu').$ Moreover let $\alpha\in Aut(A)$ satisfying 
 $\alpha(A_i)=A_i$ ($1\leq i\leq s$). Then $(A,\mu=\alpha\circ\mu',\alpha)$ is a semi-simple Hom-alternative algebra and has the unique decomposition.
 \end{enumerate}
 \end{Theorem}
 {\bf Proof.}
(1) Suppose that $(A,\mu,\alpha)$ is a finite dimensional semi-simple Hom-alternative algebra.Then $A$ has the decomposition $A=\oplus_i^s A_i$ where 
 where $A_i$($1\leq i\leq s$) are simple two-sided Hom-ideal of $(A,\mu,\alpha).$ Then $(A_i,\mu,\alpha|_{A_i})$ ($1\leq i\leq s)$ are simple finite dimensional Hom-alternative algebras. According to the proof of Proposition \ref{haat}, $\alpha|_{A_i}$
 is invertible and therefore $\alpha$ is invertible. Thus thank to Corollary \ref{chaaa}, the Hom-alternative algebra $(A,\mu,\alpha)$ is of alternative type and its induced alternative algebra is $(A,\mu')$ with 
 $\mu'=\alpha^{-1}\circ\mu.$
 
 On the other hand, by the proof of Theorem \ref{tcns1} (2), $A_i$ ($1\leq i\leq s)$ are two-sided ideal of $(A,\mu').$ Moreover $(A_i,\mu'|_{A_i})$ ($1\leq i\leq s$ are induced alternative algebra of  finite dimensional simple Hom-alternative algebras $(A_i,\mu,\alpha|_{A_i})$ ($1\leq i\leq s)$ respectively. Thank to Theorem \ref{tcns1} (1),  $(A_i,\mu')$ are semi-simple alternative algebras and can be decomposed into direct sum of isomorphic simple two-sided ideals  
 $A_i=A_{i1}\oplus A_{i2}\oplus\cdots\oplus A_{im_i}.$ It follows that 
 $(A,\mu')$ is semi-simple and has the decomposition of direct sum of simple two-sided ideals 
 $$A=A_{11}\oplus A_{12}\oplus\cdots\oplus A_{1m_1}\oplus A_{21}\oplus A_{22}\oplus\cdots\oplus A_{2m_2}\oplus\cdots\oplus A_{s1}\oplus A_{s2}\oplus\cdots A_{sm_s}.$$
(2) We know by Proposition \ref{chafaa} that $(A,\mu,\alpha)$ is a Hom-alternative algebra. Next, for all $1\leq i\leq s,$ the condition $\alpha(A_i)=A_i$ implies $$\mu(A_i,A)=\alpha(\mu'(A_i,A))\subseteq \alpha(A_i)=A_i$$ and $$\mu(A,A_i)=\alpha(\mu'(A,A_i))\subseteq \alpha(A_i)=A_i.$$ It follows that $A_i$ are two-sided Hom-ideals of $(A,\mu,\alpha).$\\
 If there exits non trivial two-sided Hom-ideal $A_{i_0}$ of $(A_i,\mu,\alpha),$ then we have
$$\mu(A_{i_0},A)=\mu(A_{i_0},A_1\oplus A_2\oplus\cdots A_s)=\mu(A_{i_0},A_i)\subseteq A_{i_0} \mbox{ and } $$
$$\mu(A,A_{i_0})=\mu(A_1\oplus A_2\oplus\cdots A_s,A_{i_0})=\mu(A_{i_0},A_i)\subseteq A_{i_0}.$$
It follows that $A_{i_0}$ is a non trivial two-sided Hom-ideal of $(A,\mu,\alpha).$ Thank to the proof of Theorem \ref{tcns1} (2), $A_{i_0}$ is also a non trivial two-sided ideal of $(A,\mu').$ Hence, $A_{i_0}$ is also a non trivial two-sided ideal of $(A_i, \mu')$ which is a contradiction. It follows that $A_i\ (i=1,\cdots, s)$ are simple two-sided Hom-ideals of $(A,\mu,\alpha)$ and therefore $(A,\mu,\alpha)$ is semi-simple and has the unique decomposition. 
\begin{Proposition}
Let $(A,\mu,\alpha)$ be a Hom-alternative algebra such that $\alpha^2=\alpha.$ Then $(A,\mu,\alpha)$ is isomorphic to the decomposition of direct sum of Hom-alternative algebras i.e.
$$A\cong (A/Ker(\alpha))\oplus Ker(\alpha).$$
\end{Proposition}
{\bf Proof.} It is clear that $(Ker(\alpha),\mu,\alpha)$ is a Hom-alternative algebra and thank to Proposition \ref{haq} since $Kerf(\alpha)$ is a two-sided Hom-ideal of $(A,\mu,\alpha),$ the quotient Hom-algebra $(A/Ker(\alpha),\bar{\mu},\bar{\alpha})$ is a Hom-alternative algebra. Now, set $A_1=(A/Ker(\alpha))\oplus Ker(\alpha)$ and define 
$\mu_1: A_1^{\times 2}\longrightarrow A_1$ and $\alpha_1: A_1\longrightarrow A_1$  by $\mu_1((\bar{x},h),(\bar{y},k)):=(\overline{\mu(x,y)},\mu(h,k))$ and $\alpha_1((\bar{x},h)):=(\overline{\alpha(x)},0).$ Then one can show that $(A_1,\mu_1,\alpha_1)$ is a Hom-alternative algebra which is a decomposition of direct sum of alternative algebras. Next, let $y\in Ker(\alpha)\cap Im(\alpha).$ Then, there exists  $x\in A$ such that $y=\alpha(x).$ Moreover, we have 
$$0=\alpha(y)=\alpha^2(x)=\alpha(x)=y.$$
It follows that $Ker(\alpha)\cap Im(\alpha)=\{0\}$ and then 
$A=Ker(\alpha)\oplus Im(\alpha)$ since for any $x\in A$ we have $x=(x-\alpha(x))+\alpha(x)$ with $x-\alpha(x)\in Ker(\alpha)$ and $\alpha(x)\in Im(\alpha).$\\
Now, let show that $(Im(\alpha),\mu,\alpha)\cong (V/Ker(\alpha),\bar{\mu},\bar{\alpha}).$ Note that, it is clear that $(Im(\alpha),\mu,\alpha)$ is a Hom-alternative algebra. Define 
$\varphi: V:Ker(\alpha)\rightarrow Im(\alpha)$ by $\varphi(\bar{x})=\alpha(x)$ for all $\bar{x}\in A/Ker(\alpha).$ Clearly, 
$\varphi$ is bijective and for all $\bar{x},\bar{y}\in A/Ker(\alpha),$ we have $$\varphi(\bar{\mu}(\bar{x},\bar{y}))=\varphi(\overline{\mu(x,u)})=\alpha(\mu(x,y))=\mu(\alpha(x),\alpha(y))=\mu(\varphi(\bar{x}),\varphi(\bar{y}))$$
$$\varphi(\bar{\alpha}(\bar{x}))=\varphi(\overline{\alpha(x)})=
\alpha^2(x)=\alpha(\varphi(\bar{x})),$$
i.e. $\varphi\circ\bar{\mu}=\bar{\mu}\circ\varphi^{\otimes 2}$ and 
$\varphi\circ\bar{\alpha}=\alpha\circ\varphi.$ It follows that 
$(Im(\alpha),\mu,\alpha)\cong (V/Ker(\alpha),\bar{\mu},\bar{\alpha})$ and therefore $V=Ker(\alpha)\oplus Im(\alpha)\cong (A/Ker(\alpha))\oplus Ker(\alpha).$
\section{Bimodule over simple Hom-alternative algebras}
In this section, we mainly study bimodules over Hom-alternative algebras of alternative type. We give a theorem about the relationship between bimodules over Hom-alternative algebras of alternative type and the ones over their induced alternative algebras. Moreover, some relevant propositions about bimodule over Hom-alternative algebras also displayed. As consequence, interesting results about bimodules over finite dimensional simple Hom-alternative algebras are given.
 \begin{Definition} Let $\mathcal{A}'=(A,\mu)$ be any algebra and $V$ be a $\mathbb{K}$-module.
\begin{enumerate}
\item A left (resp. right) structure map on $V$ is a morphism $\delta_l: A\otimes V\longrightarrow V,$
$a\otimes v\longmapsto a\cdot v$ (resp. $\delta_r: V\otimes A\longrightarrow V,$
$v\otimes a\longmapsto v\cdot a$) of Hom-modules.
\item Let $\delta_l$ and $\delta_r$ be structure maps on $V.$ Then the module associator of $V$ is a trilinear map $(,,)_{A,V}$ defined as:
\begin{eqnarray}
(,,)_{A,V}\circ Id_{V\otimes A\otimes A}=\delta_r\circ(\delta_r\otimes Id_A)-\delta_l\circ(Id_V\otimes\mu)\nonumber\\
(,,)_{A,V}\circ Id_{A\otimes V\otimes A}=\delta_r\circ(\delta_l\otimes Id_A)-\delta_l\circ(Id_A\otimes\delta_r)\nonumber\\
(,,)_{A,V}\circ Id_{A\otimes A\otimes V}=\delta_l\circ(\mu\otimes Id_V)-
\delta_l\circ(Id_A\otimes\delta_l)\nonumber
\end{eqnarray}
\end{enumerate}
\end{Definition}
A  bimodule over alternative algebras is given in \cite{Jacob2},\cite{Schaf1}.
\begin{Definition} \cite{Jacob2},\cite{Schaf1}
Let $\mathcal{A}'=(A,\mu)$ be an alternative algebra.\\
(i) An alternative $\mathcal{A}'$-bimodule is a $\mathbb{K}$-module $V$ that comes equipped with a (left) structure map  $\delta_l:A\otimes V\longrightarrow V$ ($\delta_l(a\otimes v)=a\cdot v$) and a (right) structure map $\delta_r:V\otimes A\longrightarrow V$ ($\delta_r(v\otimes a)=v\cdot a$) 
such that the following equalities:
\begin{eqnarray}
(a,v,b)_{A,V}=-(v,a,b)_{A,V}=(b,a,v)_{A,V}=-(a,b,v)_{A,V}\label{mra3}
\end{eqnarray}
\end{Definition}
 The notion of alternative bimodules  has been extended to the  Hom-alternative bimodules. More precisely, we get 
\begin{Definition} \cite{sa}
Let $\mathcal{A}=(A,\mu,\alpha_A)$ be a Hom-alternative algebra.\\
A Hom-alternative $\mathcal{A}$-bimodule is a Hom-module $(V,\alpha_V)$ that comes equipped with a (left) structure map  $\rho_l:A\otimes V\longrightarrow V$ ($\rho_l(a\otimes v)=a\cdot v$) and a (right) structure map $\rho_r:V\otimes A\longrightarrow V$ ($\rho_r(v\otimes a)=v\cdot a$) 
such that the following equalities:
\begin{eqnarray}
as_{A,V}(a,v,b)=-as_{A,V}(v,a,b)=as_{A,V}(b,a,v)=-as_{A,V}(a,b,v)\label{ra3}
\end{eqnarray}
hold for all $a,b,c\in A$ and $v\in V$ where $as_{A,V}$ is 
the module Hom-associator of the Hom-module $(V,\alpha_V)$ defined by
\begin{eqnarray}
as_{A,V}\circ Id_{V\otimes A\otimes A}=\rho_r\circ(\rho_r\otimes\alpha_A)-
\rho_l\circ(\alpha_V\otimes\mu)\nonumber\\
as_{A,V}\circ Id_{A\otimes V\otimes A}=\rho_r\circ(\rho_l\otimes\alpha_A)-
\rho_l\circ(\alpha_A\otimes\rho_r)\nonumber\\
as_{A,V}\circ Id_{A\otimes A\otimes V}=\rho_l\circ(\mu\otimes\alpha_V)-
\rho_l\circ(\alpha_A\otimes\rho_l)\nonumber
\end{eqnarray}
\end{Definition}
\begin{Remark} If $\alpha_A=Id_A$ and $\alpha_V=Id_V,$ the notion of Hom-alternative bimodule  is reduced to the one of alternative bimodule .
\end{Remark}
\begin{Theorem} \label{bHaba}
\begin{enumerate}
\item Let $\mathcal{A}=(A,\mu,\alpha_A)$ be a Hom-alternative algebra of alternative type with $\mathcal{A}'=(A,\mu')$ its induced alternative algebra and $(V,\alpha_V)$ be a Hom-alternative $\mathcal{A}$-bimodule with the structure maps $\rho_l$ and $\rho_r$ such that $\alpha_V$ is invertible. Then $V$ is an alternative $\mathcal{A}'$-bimodule with the structures maps $\delta_l=\alpha_V^{-1}\circ\rho_l$ and $\delta_r=\alpha_V^{-1}\circ\rho_r$
\item Let $\mathcal{A}'=(A,\mu')$ be an alternative algebra and  $\mathcal{A}=(A,\mu=\alpha_A\circ\mu',\alpha_A)$  a corresponding Hom-alternative algebra. Let $V$ be an alternative $\mathcal{A}'$-bimodule  with the structure maps $\delta_l$ and $\delta_r$ and $\alpha_V\in End(V)$ such that $\alpha_V\circ\delta_l=\delta_l\circ(\alpha_A\otimes \alpha_V)$
and $\alpha_V\circ\delta_r=\delta_r\circ(\alpha_V\otimes \alpha_A).$ Then $(V,\alpha_V)$ is a Hom-alternative $\mathcal{A}$-bimodule   with the structures maps $\rho_l=\alpha_V\circ\delta_l$ and $\rho_r=\alpha_V\circ\delta_r$ where $\alpha_A$ is a morphism of $(A,\mu').$
\end{enumerate}
\end{Theorem}
{\bf Proof.}
\begin{enumerate}
\item Let $(V,\alpha_V)$ be an alternative  $\mathcal{A}$-bimodule  with the structure maps $\rho_l$ and $\rho_r$ such that $\alpha_V$ is invertible. Note that the structure maps $\delta_l=\alpha_V^{-1}\circ\rho_l$ and $\delta_r=\alpha_V^{-1}\circ\rho_r$ satisfy the condition
\begin{eqnarray}
\alpha_V\circ\delta_l=\delta_l\circ(\alpha_A\otimes \alpha_V)
\mbox{ and } 
\alpha_V\circ\delta_r=\delta_r\circ(\alpha_V\otimes \alpha_A) \label{sc}
\end{eqnarray}
  Therefore
\begin{eqnarray}
as_{A,V}\circ Id_{A\otimes V\otimes A}&=&\rho_r\circ(\rho_l\otimes\alpha_A)-
\rho_l\circ(\alpha_A\otimes\rho_r)\nonumber\\
&=&\alpha_V\circ\delta_r\circ(\alpha_V\circ\delta_l\otimes\alpha_A)-
\alpha_V\circ\delta_l\circ(\alpha_A\otimes\alpha_V\circ\delta_r)\nonumber\\
&=& \alpha_V^2\circ\delta_r\circ(\delta_l\otimes Id_A)-
\alpha_V^2\circ\delta_l\circ(Id_A\otimes \delta_r) \mbox{ (by (\ref{sc} )  )}\nonumber\\
&=& \alpha_V^2\circ(,,)_{A,V}\circ Id_{A\otimes V\otimes A}\nonumber
\end{eqnarray}
Thus  $(,,)_{A,V}\circ Id_{A\otimes V\otimes A}=(\alpha_V^2)^{-1}\circ as_{A,V}\circ Id_{A\otimes V\otimes A}.$ Similarly, we get
\begin{eqnarray}
(,,)_{A,V}\circ Id_{A\otimes A\otimes V}=(\alpha_V^2)^{-1}\circ as_{A,V}\circ Id_{A\otimes A\otimes V}\nonumber\\
(,,)_{A,V}\circ Id_{V\otimes A\otimes A}=(\alpha_V^2)^{-1}\circ as_{A,V}\circ Id_{V\otimes A\otimes A}\nonumber
\end{eqnarray}
Finally, the condition (\ref{mra3}) follows from the condition (\ref{ra3}).
\item Similar  to (1).
\end{enumerate}
\begin{Corollary}
Let $\mathcal{A}=(A,\mu,\alpha_A)$ be a finite dimensional simple Hom-alternative algebra and $(V,\alpha_V)$ be a Hom-alternative $\mathcal{A}$-bimodule with the structure maps $\rho_l$ and $\rho_r$ such that $\alpha_V$ is invertible. Then $V$ is an alternative $\mathcal{A}'$-bimodule  with the structures maps $\delta_l=\alpha_V^{-1}\circ\rho_l$ and $\delta_r=\alpha_V^{-1}\circ\rho_r$ where $\mathcal{A}'=(A,\mu')$ is the induced alternative algebra.
\end{Corollary}
\begin{Definition}
Let $\mathcal{A}=(A,\mu,\alpha_A)$ be a Hom-alternative algebra and $(V,\alpha_V)$ be a Hom-alternative $\mathcal{A}$-bimodule with the structure $\rho_l$ and $\rho_r.$ 
\begin{enumerate}
\item A subspace $V_0$ of $V$ is called an $\mathcal{A}$-subbimodule of $(V,\alpha_V)$ if $\alpha_V(V_0)\subseteq V_0,$ $\rho_l(A\otimes V_0)\subseteq V_0$ and 
$\rho_r(V_0\otimes A)\subseteq V_0.$
\item  The Hom-alternative $\mathcal{A}$-module $(V,\alpha_V)$ is said to be irreducible if it has no non trivial $\mathcal{A}$-subbimodules and completely reducible if 
$V=V_1\oplus V_2\oplus\cdots\oplus V_s$ where $V_i$ ($1\leq i\leq s$) are irreducible $\mathcal{A}$-subbimodules of $(V,\alpha_V).$
\end{enumerate}
\end{Definition}
\begin{Proposition}
Let $\mathcal{A}=(A,\mu,\alpha_A)$ be a Hom-alternative algebra and $(V,\alpha_V)$ a Hom-alternative  
$\mathcal{A}$-bimodule with the structure maps $\rho_l$ and $\rho_r.$ Then 
$Ker(\alpha_V)$ is an $\mathcal{A}$-subbimodule of $(V,\alpha_V).$ 
Moreover if $\alpha_A$ is surjective then $Im(\alpha_V)$ is an $\mathcal{A}$-subbimodule of $(V,\alpha_V)$ and we have the isomorphism of $\mathcal{A}$-bimodules
 $\bar{\alpha_V}: V/Ker(\alpha_V)\rightarrow Im(\alpha_V).$
\end{Proposition}
{\bf Proof.} Obvious, we have $\alpha_V(Ker(\alpha_V))\subseteq Ker(\alpha_V).$ Next, let $(v,a)\in Ker(\alpha_V)\times A.$ Then we get
$\alpha_V(\rho_l(a\otimes v))=\rho_l(\alpha_A(a)\otimes\alpha_V(v))=0$ and 
 $\alpha_V(\rho_r(v\otimes a))=\rho_r(\alpha_V(v)\otimes\alpha_A(a))=0$ since $\alpha_V(v)=0.$ Therefore $Ker(\alpha_V)$ is an $\mathcal{A}$-subbimodule of $(V,\alpha_V).$ 
 
 Similarly, it is obvious that $\alpha_V(Im(\alpha_V))\subseteq Im(\alpha_V)).$ Let $(v,a)\in Im(\alpha_V)\times A.$
 Then there exits $v'\in V$ and if $\alpha_A$ is surjective $a'\in A$  such that $v=\alpha_V(v')$ and $a=\alpha_A(a').$ Therefore 
 $\rho_l(a\otimes v)=\rho_l(\alpha_A(a')\otimes\alpha_V(v'))=\alpha_V(\rho_l(a'\otimes v'))\in Im(\alpha_V)$ and $\rho_r(v\otimes a)=\rho_r(\alpha_V(v')\otimes\alpha_A(a'))=\alpha_V(\rho_r(v'\otimes a'))\in Im(\alpha_V).$ Thus $Im(\alpha_V)$ is an $\mathcal{A}$-subbimodule of $(V,\alpha_V).$
 
 Finally, if define the map $\bar{\alpha_V}: V/Ker(\alpha_V)\rightarrow Im(\alpha_V)$ by $\bar{\alpha_V}(\bar{v})=\alpha_V(v),$ then it is easy to prove that $\bar{\alpha_V}$ is an isomorphism.
 \begin{Corollary}
 Let $\mathcal{A}=(A,\mu,\alpha_A)$ be a Hom-alternative algebra and $(V,\alpha_V)$ be a finite dimensional irreducible Hom-alternative $\mathcal{A}$-bimodule. Then $\alpha_V$ is invertible.
 \end{Corollary}
 \begin{Proposition}
Let $\mathcal{A}=(A,\mu,\alpha_A)$ be a Hom-alternative algebra of alternative type and $(V,\alpha_V)$ a Hom-alternative
$\mathcal{A}$-bimodule with the structure maps $\rho_l$ and $\rho_r$ such that $\alpha_V$ is invertible.  If the alternative $\mathcal{A}'$-bimodule $V$ over the induced alternative algebra $\mathcal{A}'=(A,\mu')$ with the structures maps $\delta_l=\alpha_V^{-1}\circ\rho_l$ and $\delta_r=\alpha_V^{-1}\circ\rho_r$ is irreducible, then the Hom-alternative $\mathcal{A}$-bimodule $(V,\alpha_V)$ is also irreducible.
\end{Proposition}
{\bf Proof.} Assume that the Hom-alternative $\mathcal{A}$-bimodule $(V,\alpha_V)$ is reducible. Then there exists a non trivial subspace $V_0$ such that $(V_0, \alpha_V|_{V_0})$ is an $\mathcal{A}$-subbimodule of $(V,\alpha_V).$  Therefore $\alpha_V(V_0)\subseteq V_0,$ 
$\rho_l(a\otimes v))\in V_0$ and $\rho_r(v\otimes a)\in V_0$ for all 
$(a,v)\in A\times V_0.$  Hence  
$\delta_l(a\otimes v))\in \alpha_V^{-1}(V_0)=V_0$ and $\delta_r(v\otimes a)\in \alpha_V^{-1}(V_0)=V_0.$  Thus  $V_0$ is a non trivial $\mathcal{A}'$-subbimodule of $V,$ contradiction. It follows that $(V,\alpha_V)$ is an irreducible Hom-alternative $\mathcal{A}$-bimodule.
\begin{Corollary}
Let $\mathcal{A}=(A,\mu,\alpha_A)$ be a finite dimensional simple Hom-alternative algebra and $(V,\alpha_V)$ a Hom-alternative $\mathcal{A}$-bimodule with the structure maps $\rho_l$ and $\rho_r$ such that $\alpha_V$ is invertible.  If the alternative $\mathcal{A}'$-bimodule $V$ over the induced alternative algebra $\mathcal{A}'=(A,\mu')$ with the structures maps $\delta_l=\alpha_V^{-1}\circ\rho_l$ and $\delta_r=\alpha_V^{-1}\circ\rho_r$ is irreducible, then the Hom-alternative $\mathcal{A}$-bimodule $(V,\alpha_V)$ is also irreducible.
\end{Corollary}
Let recall the following result from \cite{Schaf1} which is very useful for the next result.
\begin{Theorem}\cite{Schaf1}
Let $\mathcal{A}'=(A,\mu')$ be a semi-simple alternative algebra. Then any  representation $(S,T)$ of $\mathcal{A}'$ is completely reducible. 
\end{Theorem}
Since the notion of completely reducible representation of $\mathcal{A}',$ is equivalent to the notion of 
a completely reducible alternative $\mathcal{A}'$-bimodule,  thank to  Theorem \ref{ts} (1) and Theorem \ref{bHaba} (1), we get the following important result.
\begin{Corollary}
Let $\mathcal{A}=(A,\mu,\alpha)$ be a finite dimensional semi-simple Hom-alternative algebra. Then any Hom-alternative $\mathcal{A}$-bimodule $(V,\alpha_V)$ such that $\alpha_V$ is invertible, is is completely reducible.
\end{Corollary}

Author's addresses\\
Sylvain Attan, syltane2010@yahoo.fr\\
D\'epartement de Math\'ematiques,
Universit\'e d'Abomey-Calavi,
01 BP 4521, Cotonou 01, B\'enin.
\end{document}